\documentclass[12pt,twoside]{article}
\usepackage{mathrsfs}

\usepackage{amssymb,amsmath,amsthm,epsfig,latexsym}
\usepackage[numbers,sort&compress]{natbib}

\usepackage{txfonts}
\usepackage{amsfonts}
\usepackage{color}

\usepackage{geometry}

\usepackage [latin1]{inputenc}
\newtheorem{theorem}{Theorem}[section]

\newtheorem{lem}{Lemma}[section]

\numberwithin{equation}{section}

\begin{document}
\renewcommand{\baselinestretch}{1.15}
\begin{center}
{\Large\bf A class of variable-order fractional $p(\cdot)$-Kirchhoff type systems}
\vskip.2in {\small Yong Wu$^{1}$  \qquad Zhenhua Qiao$^{2}$ \qquad  Mohamed Karim Hamdani$^{3,4,5}$\qquad  Bingyu Kou$^{6}$\footnote{Corresponding author: koubei@163.com (B.Y. Kou)}\qquad  Libo Yang$^{7}$}\\[1mm]
{\scriptsize $^{1}$School of Tourism Date, Guilin Tourism University, Guilin, 541006, China.}\\
{\scriptsize $^{2}$School of Electronic and Information Engineering, Jiangxi Industry Polytechnic College, Nanchang 330095, China.}\\
{\scriptsize $^{3}$Science and technology for defense Laboratory  LR19DN01, Military Research Center, Aouina, Tunisia.}
\\{\scriptsize $^{4}$Military School of Aeronautical Specialities, Sfax, Tunisia.}\\
{\scriptsize $^{5}$Mathematics Department, University of Sfax, Faculty of Science of Sfax, Sfax, Tunisia.}\\
{\scriptsize $^{6}$Department of Basic Courses, The Army Engineering University of PLA, Nanjing, 211101, China.}\\
{\scriptsize $^{7}$Faculty of Mathematics and Physics, Huaiyin Institute of Technology, Huai'an 223003, China.}
\end{center}

\vspace{0.05in}

{\footnotesize \noindent {\bf Abstract}~   This paper is concerned with an elliptic system of Kirchhoff type, driven by the variable-order fractional $p(x)$-operator. With the help of the direct variational method and Ekeland variational principle, we show the existence of a weak solution.  This is our first attempt to study this kind of system, in the case of variable-order fractional variable exponents. Our main theorem extends in several directions previous results.\\

}

\vskip.2in

\section{Introduction}\label{sec1}
In this article, we discuss the following variable-order fractional $p(\cdot)$-Kirchhoff type system
\begin{equation}\label{1.1}
\left\{\begin{array}{cl}
&M_{1}\left(\displaystyle\iint_{\mathbb R^{2N}}\frac{1}{p(x,y)}\displaystyle{\frac{|u(x)-u(y)|^{p(x,y)}}{|x-y|^{N+p(x,y)s(x,y)}}}dxdy\right)(-\Delta)^{s(\cdot)}_{p(\cdot)}u(x)
=f(u,v)+a(x)  ~~~~ {\rm in}~~\Omega,\\
&M_{2}\left(\displaystyle\iint_{\mathbb R^{2N}}\frac{1}{p(x,y)}\displaystyle{\frac{|v(x)-v(y)|^{p(x,y)}}{|x-y|^{N+p(x,y)s(x,y)}}}dxdy\right)(-\Delta)^{s(\cdot)}_{p(\cdot)}v(x)
=g(u,v)+b(x)  ~~~~~ {\rm in}~~\Omega,\\
&u=v=0 ~~~~~~~~~~~~~~~~~~~~~~~~~~~~~~~~~~~~~~~~~~~~~~~~~~~~~~~~~~~~~~~~~~~~~~~~~~~~~~~~~~~~~~~~~~~~~~~~~~~{\rm in}~~\mathbb{R}^{N}\setminus \Omega,
\end{array}\right.
\end{equation}
where $\Omega\subset\mathbb{R}^{N}$ is a bounded smooth domain with $N>p(x,y)s(x,y)$ for any $(x,y)\in\overline{\Omega}\times\overline{\Omega}$.
Here, the main operator $(-\Delta)^{s(\cdot)}_{p(\cdot)}$ is the variable-order fractional $p(\cdot)$-Laplacian given by
\begin{equation}\label{mio}
(-\Delta)^{s(\cdot)}_{p(\cdot)}\varphi(x)=P.V.\int_{\mathbb{R}^{N}}\frac{|\varphi(x)-\varphi(y)|^{p(x,y)-2}(\varphi(x)-\varphi(y))}{|x-y|^{N+p(x,y)s(x,y)}}dy, ~~x\in\mathbb{R}^{N},
\end{equation}
along any $\varphi\in C_0^\infty(\mathbb R^N)$, where P.V. denotes the Cauchy principle value.

From now on, in order to simplify the notation, we denote
$$
\begin{gathered}
 \quad p^{-}=\min\limits_{(x,y)\in\mathbb R^{2N}}p(x,y),\quad p^{+}=\max\limits_{(x,y)\in\mathbb R^{2N}}p(x,y), \quad s^{-}=\min\limits_{(x,y)\in\mathbb R^{2N}}s(x,y),\quad s^{+}=\max\limits_{(x,y)\in\mathbb R^{2N}}s(x,y).
\end{gathered}
$$
We will assume that $M_{1}, M_{2}:\mathbb R^+\to\mathbb R^+$ are  \textit{continuous} functions satisfying the condition\\
$(M):$ \textit{there exist $m>0$ and $\gamma>\frac{1}{p^{-}}$ such that}
$$M_{1}(t),M_{2}(t)>mt^{\gamma-1}, ~\textit{for all} ~t>0.$$
Note that the Kirchhoff functions $M_{1}, M_{2}$ may be singular at $t=0$ for $\gamma\in(0,1)$;

Moreover, $H:\mathbb R^{2}\to\mathbb R$ is a $C^{1}$-function verifying\\
$(Z_{1}):$
$$ \frac{\partial H}{\partial u}(u,v)=f(u,v)~~and ~~~\frac{\partial H}{\partial v}(u,v)=g(u,v)~~for~all~~(u,v)\in \mathbb{R}^{2};$$
$(Z_{2}):$
\textit{there exists $K>0$ such that}
$$H(u,v)=H(u+K,v+K)~~for~all~~(u,v)\in \mathbb{R}^{2};$$
Finally, we suppose that\\
\noindent $(AB):$ \textit{$a(x),b(x)\in L^{q(x)}(\Omega), \frac{1}{p(x)}+ \frac{1}{q(x)}=1, 1<q(x)<p^{*}_{s}(x)$.}
where $p^{*}_{s}(x)=N\overline{p}(x)/(N-\overline{s}(x)\overline{p}(x))$,  $\overline{p}(x)=p(x,x)$, $\overline{s}(x)=s(x,x)$\\
\noindent $(PS):$ $p(\cdot):\mathbb R^{2N}\to(1,\infty)$ is a \textit{continuous} function fulfilling $0<s^{-}<s^{+}<1<p^{-}\leq p^{+}$,
and $p(\cdot)$ and $s(\cdot)$ are symmetric, that is,  $p(x,y)=p(y,x)$ and  $s(x,y)=s(y,x)$ for any $(x,y)\in\mathbb R^{2N}$.

On the one hand, when $s(\cdot)\equiv 1$, the operators in \eqref{1.1} reduce to the integer order, i.e., $p(\cdot)$-Laplacian $\Delta_{p(\cdot)}$. This kind of variable exponents problem has a wide range of real applications, such as electrorheological fluids (see \cite{eg}), elastic mechanics (\cite{vvz}), image restoration (\cite{ycslmr}) and so on. For this kind of operator combined with kirchhoff function system problem,  we recall \cite{1,2,3,4,5,HR}, for example, S. Boulaaras (\cite{1}) at al. studied the existence of positive solutions of a p(x)-Kirchhoff system by using sub-super solutions concepts. A very interesting question arose is that whether there are other ways to solve this class of problems. And also, can we consider the non-local variable-order case? We known that the fractional variable order derivatives proposed by Lorenzo and
Hartley in \cite{X23} appeared in different nonlinear diffusion processes. Subsequently, many results of the variable order problem have appeared in the literature \cite{nuovo,lz,zzz}.

Of course, when $p(x)\equiv p (\text{or}~ p=2)$ and $s(\cdot)\equiv s$, the operators in \eqref{1.1} reduce to the classical non-local fractional $p$-Laplacian i.e., $(-\Delta)^{s}_{p}$ ($(-\Delta)^{s}$). Similarly, for Kirchhoff type system cases, the papers \cite{6,7,8,9,10,11,12} introduces a lot of related work in recent years, where many authors studied existence and multiplicity of solutions by applying variational methods. For instance, based on the three critical points theorem, E. Azroul at al. (\cite{7}) discussed an elliptic system  with homogeneous Dirichlet boundary conditions, they obtained the existence of three weak solutions.

On the other hand, it is worth mentioning that Kirchhoff in 1883 (see\cite{kkk}) presented a stationary verion of differential equation, the so-called kirchhoff equation
$$\rho\frac{\partial^{2} v}{\partial t^{2}}-\Big(\frac{p_{0}}{l}+\frac{e}{2L}\int_{0}^{L}|\frac{\partial v}{\partial x}|^{2}dx\Big)\frac{\partial^{2}v}{\partial x^{2}}=0,$$
where $\rho,l, e,L, p_{0}$ are positive constants which represent the corresponding physical meanings. It's a generalization of D'Alembert equation.
It's very interesting to combine this model with various operators due to its nonlocal nature.

Inspired by the above works, we consider a new fractional Kirchhoff type system \eqref{1.1}. As far as we know, this is the first attempt on variable-order fractional situations to study a bi-non-local problem with variable exponent. In order to overcome the difficulty, we use the direct variational method and Ekeland variational principle to deal with it. Our result is new to the variable-order fractional system with variable exponent.

Now, we give the main result of this paper, our energy functional $\mathcal{I}$ will be introduced in Section 2.
\begin{theorem}\label{T1.1}
Let $\Omega$ be a bounded smooth domain of $\mathbb{R}^{N}$, with $N>p(x,y)s(x,y)$ for any $(x,y)\in\overline{\Omega}\times\overline{\Omega}$, where $p(\cdot)$ and $s(\cdot)$ verify $(PS)$.
Assume that $(M)$, $(AB)$ and $(Z_{1})-(Z_{2})$ are satisfied.
Then, problem \eqref{1.1} admits a weak solution if $\mathcal{I}$ is differentiable at $(u_{0}, v_{0})$.
\end{theorem}

The paper is organized as follows. In Section \ref{sec2}, we state some interesting properties of variable exponent Lebesgue spaces and variable-order fractional Sobolev spaces with variable exponent. In Section \ref{sec3}, we prove the functional $\mathcal{I}$ is bounded from below and give the proof of Theorem \ref{T1.1}.

\section{Abstract framework}\label{sec2}
In this section, first of all, we recall some basic properties about the variable exponent Lebesgue spaces in \cite{RD} and variable-order fractional Sobolev spaces. Secondly, we give some necessary lemmas that will be used in this paper. Finally, we introduce the definition of weak solutions for problem \eqref{1.1} and build the corresponding energy functional.
Consider the set $$C_{+}(\overline{\Omega})=\bigg\{p\in C(\overline{\Omega}):p(x)>1 ~~\text{for~all}~~ x\in \overline{\Omega}\bigg\}.$$
For any $p\in C_{+}(\overline{\Omega}), $ we define the variable exponent Lebesgue space as

$$L^{p(\cdot)} (\Omega)=\Bigg\{u: ~\text{the function}~ u: \Omega\rightarrow \mathbb{R}~~ \text{is measurable}, ~~\int_{\Omega}|u(x)|^{p(x)}dx<\infty\Bigg\},$$
the vector space endowed with the $Luxemburg norm$,
$$\|u\|_{p(\cdot)}=\inf\Bigg\{\lambda>0:\int_{\Omega}\Big|\frac{u(x)}{\lambda}\Big|^{p(x)}dx\leq1\Bigg\}.$$
Then $(L^{p(\cdot)}(\Omega),\|\cdot\|_{p(\cdot)})$ is a separable reflexive Banach space, see \cite[Theorem 2.5]{39}.
Let $q\in C_{+}(\overline{\Omega})$ be the conjugate exponent of $p$, that is
$$\frac{1}{p(x)}+\frac{1}{q(x)}=1,\qquad\mbox{for all }x\in\overline{\Omega}.$$
Then we have the following H\"{o}lder inequality, whose proof can be found in \cite[Theorem 2.1]{39}.

\begin{lem}\label{L2.1}
Assume that $u\in L^{p(\cdot)}(\Omega)$ and $v\in L^{q(\cdot)}(\Omega)$, then
\begin{equation*}
\Big|\int_{\Omega}uvdx\Big|\leq\Big(\frac{1}{p^{-}}+\frac{1}{q^{-}}\Big)\|u\|_{p(\cdot)}\|v\|_{q(\cdot)}\leq2\|u\|_{p(\cdot)}\|v\|_{q(\cdot)}.
\end{equation*}
\end{lem}

The variable-order fractional Sobolev spaces with variable exponent is defined by
$$W^{s(\cdot),p(\cdot)}(\Omega)=\Bigg\{u\in L^{\overline{p}(\cdot)}(\Omega): \iint_{\Omega\times\Omega}\frac{|u(x)-u(y)|^{p(x,y)}}{|x-y|^{N+p(x,y)s(x,y)}}dxdy<\infty\Bigg\}$$
with the norm $\|u\|_{s,p(\cdot)}=\|u\|_{p(\cdot)}+[u]_{s(\cdot),p(\cdot)}$, where
$$[u]_{s(\cdot),p(\cdot)}=\inf\Bigg\{\lambda>0:\iint_{\Omega\times\Omega}\frac{|u(x)-u(y)|^{p(x,y)}}{\lambda^{p(x,y)}|x-y|^{N+p(x,y)s(x,y)}}dxdy<1\Bigg\}.$$
For a more detailed introduction of this space we can refer to \cite{add7}.
For the reader's convenience, we now list some of the results in reference \cite{add7} which will be used in our paper.
We define the new variable order fractional Sobolev spaces with variable exponent
$$X=\Bigg\{u:\mathbb{R}^{N}\rightarrow\mathbb{R}:~~ u|_{\Omega}\in L^{\overline{p}(\cdot)}(\Omega), \iint_{Q}\frac{|u(x)-u(y)|^{p(x,y)}}{\lambda^{p(x,y)}|x-y|^{N+p(x,y)s(x,y)}}dxdy<\infty,~~ \text{for some} ~~\lambda>0\Bigg\},$$
where $Q:=\mathbb{R}^{2N}\backslash(\Omega^{c}\times\Omega^{c})$. The space $X$ is endowed with the norm
$$\|u\|_{X}=\|u\|_{p(\cdot)}+[u]_{X},$$
where $$[u]_{X}=\inf\Bigg\{\lambda>0:\iint_{Q}\frac{|u(x)-u(y)|^{p(x,y)}}{\lambda^{p(x,y)}|x-y|^{N+p(x,y)s(x,y)}}dxdy<1\Bigg\}.$$
We know that the norms $\|\cdot\|_{s,p(\cdot)}$ and $\|\cdot\|_{X}$ are not the same due to the fact that $\Omega\times\Omega\subset Q$  and $\Omega\times\Omega\neq Q$. This makes the variable order fractional Sobolev space with variable exponent $W^{s(\cdot),p(\cdot)}(\Omega)\times W^{s(\cdot),p(\cdot)}(\Omega)$ not sufficient for investigating the class of problems like \eqref{1.1}.

For this, we set space as
$$X_{0}=\Big\{u\in X:~ u=0 ~~a.e.~~ \text{in}~~ \mathbb{R}^{N}\backslash\Omega\Big\}.$$
The space $X_{0}$ is a separable reflexive Banach space, see \cite{m5}, with respect to the norm
\begin{align*}
\|u\|_{X_{0}}&=\inf\Bigg\{\lambda>0: \iint_{Q}\frac{|u(x)-u(y)|^{p(x,y)}}{\lambda^{p(x,y)}|x-y|^{N+p(x,y)s(x,y)}}dxdy
=\iint_{\mathbb R^{2N}}\frac{|u(x)-u(y)|^{p(x,y)}}{\lambda^{p(x,y)}|x-y|^{N+p(x,y)s(x,y)}}dxdy<1\Bigg\},
\end{align*}
where last equality is a consequence of the fact that $u=0$ a.e. in $\mathbb{R}^{N}\backslash\Omega$.\\

In the following Lemma, we give a compact embedding result. For the proof we refer the reader to \cite{add7}.

\begin{lem}\label{L2.2}
Let  $\Omega \subset \mathbb{R}^n$ be a smooth bounded domain and $s(\cdot)\in (0,1)$. Let $p(x,y)$ be continuous variable exponents with $s(x,y)p(x,y)<N$ for $(x,y)\in \overline{\Omega}\times\overline{\Omega}$. Assume that $q: \overline{\Omega}\longrightarrow(1,\infty)$ is a continuous function such that
$$p^*_{s}(x) > q(x) \geq q^->1,\;\;\mbox{for~all}\;\; x\in \overline{\Omega}.$$
Then, there exists a constant $C=C(N,s,p,q,\Omega)$ such that for every $u\in X_{0}$, it holds that
$$\|u\|_{q(x)}\leq C\|u\|_{X_{0}}.$$
The space $X_{0}$ is continuously embedded in $L^{q(x)}(\Omega)$. Moreover, this embedding is compact.
\end{lem}

We define the fractional modular function $\varrho_{p(\cdot)}^{s(\cdot)}:X_{0}\rightarrow\mathbb{R}$, by
$$\varrho_{p(\cdot)}^{s(\cdot)}(u)=\iint_{\mathbb R^{2N}}\frac{|u(x)-u(y)|^{p(x,y)}}{|x-y|^{N+p(x,y)s(x,y)}}dxdy.$$

We also have the next result of \cite[Proposition 2.2]{add7}.
\begin{lem}\label{L2.3}
Assume that $u\in X_0$ and $\{u_j\}\subset X_0$, then
\begin{align*}
\noindent(1)~~&\|u\|_{X_{0}}<1(resp.=1,>1)\Leftrightarrow \varrho_{p(\cdot)}^{s(\cdot)}(u)<1(resp.=1,>1),\\
\noindent(2)~~&\|u\|_{X_{0}}<1\Rightarrow\|u\|_{X_{0}}^{p^{+}}\leq\varrho_{p(\cdot)}^{s(\cdot)}(u)\leq\|u\|_{X_{0}}^{p^{-}},\\
\noindent(3)~~&\|u\|_{X_{0}}>1\Rightarrow\|u\|_{X_{0}}^{p^{-}}\leq\varrho_{p(\cdot)}^{s(\cdot)}(u)\leq\|u\|_{X_{0}}^{p^{+}},\\
\noindent(4)~~&\lim\limits_{j\rightarrow\infty}\|u_{j}\|_{X_{0}}=0(\infty)\Leftrightarrow\lim\limits_{j\rightarrow\infty}\varrho_{p(\cdot)}^{s(\cdot)}(u_{j})=0(\infty),\\
\noindent(5)~~&\lim\limits_{j\rightarrow\infty}\|u_{j}-u\|_{X_{0}}=0\Leftrightarrow\lim\limits_{j\rightarrow\infty}\varrho_{p(\cdot)}^{s(\cdot)}(u_{j}-u)=0.
\end{align*}
\end{lem}
Finally, we define our workspace $S=X_{0}\times X_{0}$.which is endowed with the norm
\begin{align*}
\|(u,v)\|_{S}=\|u\|_{X_{0}}+\|v\|_{X_{0}}.
\end{align*}
We say that a pair of functions $(u,v)\in S$ is the weak solution of problem \eqref{1.1}, if for all $(\phi,\varphi)\in S$ one has
\begin{align*}
&M_{1}\Big(\delta_{p(\cdot)}(u)\Big)\iint_{\mathbb R^{2N}}\!\frac{|u(x)-u(y)|^{p(x,y)-2}(v(x)-v(y))(\phi(x)-\phi(y))}{|x-y|^{N+p(x,y)s(x,y)}}dxdy=\int_{\Omega}\left((f(u,v)+a(x)\right)\phi dx,\nonumber\\
&M_{2}\Big(\delta_{p(\cdot)}(v)\Big)\iint_{\mathbb R^{2N}}\!\frac{|v(x)-v(y)|^{p(x,y)-2}(v(x)-v(y))(\varphi(x)-\varphi(y))}{|x-y|^{N+p(x,y)s(x,y)}}dxdy=\int_{\Omega}\left((g(u,v)+b(x)\right)\varphi dx,
\end{align*}
 where
$$\delta_{p(\cdot)}(u)=\iint_{\mathbb R^{2N}}\frac{1}{p(x,y)}\frac{|u(x)-u(y)|^{p(x,y)}}{|x-y|^{N+p(x,y)s(x,y)}}dxdy.$$
Let us consider the following functional associated to problem \eqref{1.1}, defined by $\mathcal{I}:S\rightarrow\mathbb{R}$
\begin{align*}
\mathcal{I}(u,v)=\widetilde{M_{1}}\Big(\delta_{p(\cdot)}(u)\Big)-\widetilde{M_{2}}\Big(\delta_{p(\cdot)}(v)\Big)-\int_{\Omega}H(u,v)dx-\int_{\Omega}a(x)udx-\int_{\Omega}b(x)vdx,
\end{align*}
for all $(u,v)\in S$, where $\widetilde{M_{i}}(t)=\int_{0}^{t}M_{i}(\tau)d\tau.$
Obviously, the continuity of $M$ yields that $\mathcal{I}$ is well defined and of class $C^{1}$ on $S\setminus\{0,0\}$. Furthermore, for every $(u,v)\in S\setminus\{0,0\},$ the derivative of $\mathcal{I}$ is given by
\begin{align*}
\langle \mathcal{I}^{'}(u,v),(\phi,\varphi) \rangle=& M_{1}\Big(\delta_{p(\cdot)}(u)\Big)\iint_{\mathbb R^{2N}}\frac{|u(x)-u(y)|^{p(x,y)-2}(u(x)-u(y))(\phi(x)-\phi(y))}{|x-y|^{N+p(x,y)s(x,y)}}dxdy\nonumber\\&+M_{2}\Big(\delta_{p(\cdot)}(v)\Big)\iint_{\mathbb R^{2N}}\frac{|v(x)-v(y)|^{p(x,y)-2}(v(x)-v(y))(\varphi(x)-\varphi(y))}{|x-y|^{N+p(x,y)s(x,y)}}dxdy\nonumber\\
&-\int_{\Omega}\left((f(u,v)+a(x)\right)\phi dx-\int_{\Omega}\left((g(u,v)+b(x)\right)\varphi dx,
\end{align*}
for any $(\phi,\varphi)\in S$ . Therefore, the weak solution $(u,v)\in S\setminus\{0,0\}$ of problem \eqref{1.1} is a nontrivial critical point of $\mathcal{I}$.

Now, we recall the following well-known Ekeland variational principle found in \cite{4}, which will be used to prove our  conclusion, that is Theorem \ref{T1.1}.

\begin{theorem}\label{T2.1}
Let $X$ be a Banach space and $\mathcal{I}: X\rightarrow \mathbb{R}$ be a $C^{1}$ function which is bounded from below. Then, for any $\varepsilon>0$, there exists $\varpi_{\varepsilon}\in X$ such that
$$\mathcal{I}(\varpi_{\varepsilon})\leq\inf\limits_{X}\mathcal{I}+\varepsilon~~and~~\|\mathcal{I}^{'}(\varpi_{\varepsilon})\|_{X^{*}}\leq\varepsilon$$
\end{theorem}

Throughout the paper, for simplicity, we use $\{c_{i},~ i\in\mathbb N\}$ to denote different non-negative or positive constant.

\section{The main result}\label{sec3}

\begin{lem}\label{lemma1}
Under the same assumptions of Theorem \ref{T1.1}, then $\mathcal{I}$ is coercive and bounded from below.
\end{lem}
{\bf Proof.} Firstly, we know that functional $\mathcal{I}$ is well defined. Indeed, it is sufficient to prove that the functional $T:S\rightarrow\mathbb{R}$, $T(u,v)=\int_{\Omega}H(u,v)dx$, is well defined. Since $H$ is continuous on $[0,K]\times[0,K]$ and $H(u,v)=H(u+K,v+K)$ for all $(u,v)\in \mathbb{R}^{2}$, we get $|H(u,v)|\leq c_{1}$ for all $(u,v)\in \mathbb{R}^{2}$. Thus,
$$|T(u,v)|\leq\int_{\Omega}|H(u,v)|dx\leq c_{1}|\Omega|, ~~for~ all~~ (u,v)\in \mathbb{R}^{2},$$
i.e., $T$ is well defined, where $|\Omega|$ is the Lebesgue measure of $\Omega$. Next, we will prove that $\mathcal{I}$ is coercive and bounded from below. Let $(u,v)\in S$, we have

\begin{align*}
\mathcal{I}(u,v)=&\widetilde{M_{1}}\Big(\delta_{p(\cdot)}(u)\Big)-\widetilde{M_{2}}\Big(\delta_{p(\cdot)}(v)\Big)-\int_{\Omega}H(u,v)dx-\int_{\Omega}a(x)udx-\int_{\Omega}b(x)vdx\\
&\geq\widetilde{M_{1}}\Big(\delta_{p(\cdot)}(u)\Big)-\widetilde{M_{2}}\Big(\delta_{p(\cdot)}(v)\Big)-c_{1}|\Omega|-\int_{\Omega}a(x)udx-\int_{\Omega}b(x)vdx.
\end{align*}
By the condition $(AB)$ and Lemma \ref{L2.1}, we get
\begin{align*}
\mathcal{I}(u,v)\geq\widetilde{M_{1}}\Big(\delta_{p(\cdot)}(u)\Big)-\widetilde{M_{2}}\Big(\delta_{p(\cdot)}(v)\Big)-c_{1}|\Omega|-2\|a(x)\|_{q(x)}\|u\|_{p(x)}-2\|b(x)\|_{q(x)}\|v\|_{p(x)}.
\end{align*}
It follows from $(M)$ and Lemmas \ref{L2.2}-\ref{L2.3} that
\begin{align}\label{1}
\mathcal{I}(u,v)\geq& m\int_{0}^{\frac{1}{p^{+}}\varrho_{p(\cdot)}^{s}(u)}\tau^{\gamma-1} d\tau+ m\int_{0}^{\frac{1}{p^{+}}\varrho_{p(\cdot)}^{s(\cdot)}(v)}\tau^{\gamma-1} d\tau-c_{3}\|u\|_{X_{0}}-c_{4}\|v\|_{X_{0}}-c_{2}\nonumber\\
&=\frac{m}{\gamma(p^{+})^{\gamma}}\left((\varrho_{p(\cdot)}^{s(\cdot)}(u))^{\gamma}+(\varrho_{p(\cdot)}^{s(\cdot)}(v))^{\gamma}\right)-c_{3}\|u\|_{X_{0}}-c_{4}\|v\|_{X_{0}}-c_{2}\nonumber\\
&\geq\frac{m}{\gamma(p^{+})^{\gamma}}\left(\min\{\|u\|_{X_{0}}^{\gamma p^{-}}, \|u\|_{X_{0}}^{\gamma p^{+}}\}+\min\{\|v\|_{X_{0}}^{\gamma p^{-}}, \|v\|_{X_{0}}^{\gamma p^{+}}\} \right)-\max\{c_{3},c_{4}\}(\|u\|_{X_{0}}+\|v\|_{X_{0}})-c_{2},
\end{align}
since $\gamma p^{+}>\gamma p^{-}>1$, when $\|(u,v)\|_{S}\rightarrow+\infty$, at least one of $\|u\|_{X_{0}}$ and $\|v\|_{X_{0}}$ converges to infinity. So, $\mathcal{I}$ is coercive and bounded from below. The proof of Lemma \ref{lemma1} is complete. \qed
\subsection*{ \bf  Proof of Theorem \ref{T1.1}.}
Obviously, since $\mathcal{I}\in C^{1}(S, \mathbb{R})$ is weakly lower semi-continuous and bounded from below, by means of Ekeland variational principle we have $(u_{j}, v_{j})\subset S$ such that
\begin{align}\label{2}\mathcal{I}(u_{j}, v_{j})\rightarrow\inf\limits_{S}\mathcal{I}~~ and ~~\mathcal{I}^{'}(u_{j}, v_{j})\rightarrow0.\end{align}
Furthermore, by the above expression, we get $|\mathcal{I}(u_{j}, v_{j})|\leq c_{5}$. Thus, it follows from \eqref{1} that
$$c_{6}\leq|\mathcal{I}(u_{j}, v_{j})|\leq c_{5}$$
which implies that the sequences $\{u_{j}\}$ and $\{v_{j}\}$ are bounded in $X_{0}$. So, without loss of generality, there exist subsequences  $\{u_{j}\}$ and $\{v_{j}\}$ such that $u_{j}\rightharpoonup u_{0}$ and $v_{j}\rightharpoonup v_{0}$ in $X_{0}$, and thus,
$$\int_{\Omega}a(x)u_{j}dx\rightarrow\int_{\Omega}a(x)u_{0}dx~~and~~\int_{\Omega}b(x)v_{j}dx\rightarrow\int_{\Omega}a(x)v_{0}dx.$$
According to compact embedding theorem, which is lemma \ref{L2.2}, we obtain
$$u_{j}(x)\rightarrow u_{0}(x)~~and~~v_{j}(x)\rightarrow v_{0}(x)~~a.e. ~~x\in \Omega.$$
Again,by continuity of $H$, we get
$$H(u_{j}(x), v_{j}(x))\rightarrow H(u_{0}(x), v_{0}(x))~~a.e. ~~x\in \Omega.$$
And because $H$ is bounded, we get the following convergence from the  Lebesgue dominated convergence theorem,
$$\int_{\Omega}H(u_{j}, v_{j})dx\rightarrow\int_{\Omega}H(u_{0}, v_{0})dx.$$
By \eqref{2}, we note that
\begin{align*}
\inf\limits_{S}\mathcal{I}=&\lim\mathcal{I}(u_{j}, v_{j})\\=&\lim\left(\widetilde{M_{1}}\Big(\delta_{p(\cdot)}(u_{j})\Big)-\widetilde{M_{2}}\Big(\delta_{p(\cdot)}(v_{j})\Big)-\int_{\Omega}H(u_{j},v_{j})dx-\int_{\Omega}a(x)u_{j}dx-\int_{\Omega}b(x)v_{j}dx.\right)
\end{align*}
In view of Fatou's Lemma, we have
$$\delta_{p(\cdot)}(u_{0})\leq\liminf\delta_{p(\cdot)}(u_{j})~~and~~\delta_{p(\cdot)}(v_{0})\leq\liminf\delta_{p(\cdot)}(v_{j}).$$
By the continuous monotone increasing property of $\widetilde{M_{1}}$ and  $\widetilde{M_{2}}$, we get
$$\widetilde{M_{1}}\Big(\delta_{p(\cdot)}(u_{0})\Big)\leq\lim\widetilde{M_{1}}\Big(\delta_{p(\cdot)}(u_{j})\Big)~~and~~\widetilde{M_{2}}\Big(\delta_{p(\cdot)}(v_{0})\Big)\leq\lim \widetilde{M_{2}}\Big(\delta_{p(\cdot)}(v_{j})\Big).$$
In conclusion,
\begin{align*}
\inf\limits_{S}\mathcal{I}\geq\widetilde{M_{1}}\Big(\delta_{p(\cdot)}(u_{0})\Big)-\widetilde{M_{2}}\Big(\delta_{p(\cdot)}(v_{0})\Big)-\int_{\Omega}H(u_{0},v_{0})dx-\int_{\Omega}a(x)u_{0}dx-\int_{\Omega}b(x)v_{0}dx=\mathcal{I}(u_{0},v_{0}),
\end{align*}
which implies $\mathcal{I}(u_{0},v_{0})=\inf\limits_{S}\mathcal{I}$. Thus, $(u_{0},v_{0})\in S$ is a weak solution of problem \eqref{1.1} if $\mathcal{I}$ is differentiable at $(u_{0},v_{0})$. The proof is complete.

\section*{Acknowledgment}
The third author would like to express his deepest gratitude to the Military School of Aeronautical Specialities, Sfax (ESA) for providing an excellent atmosphere for work. The fourth author is supported by the Fundamental
Research Funds for Youth Development of The Army Engineering University of PLA. (Grant No.
KYJBJQZL2003). The fifth author is supported by the Natural Science Foundation of Huaiyin Institute of Technology.( Grant/Award
Number: 20HGZ002).

\section*{Competing ~interests}
The authors declare that they have no competing interests.

\end{document}